\documentclass[11pt]{amsart}
\usepackage{amsmath,amsthm,amscd,amssymb}

\newcommand{\C}{{\mathbb C}}

\newcommand{\R}{{\mathbb R}}

\newcommand{\Q}{{\mathbb Q}}
\newcommand{\D}{{\mathbb D}}

\newcommand{\cC}{{\mathcal C}}

\newcommand{\ds}{\displaystyle}

\theoremstyle{plain}
\newtheorem{lemma}{Lemma}[section]
\newtheorem{thm}[lemma]{Theorem}
\newtheorem{prop}[lemma]{Proposition}
\newtheorem{cor}[lemma]{Corollary}

\theoremstyle{definition}

\begin{document}
\title[Smooth Siegel discs without number theory]{Smooth Siegel discs
  without number theory: A remark on a proof by Buff and Ch{\'e}ritat}
\author{Lukas Geyer}
\address{Lukas Geyer\\Montana State University\\Department of
Mathematics\\ P.O.~Box 172400\\Bozeman, MT 59717--2400\\ USA}
\email{geyer@math.montana.edu}
\date{}
\thanks{The author was supported by a Feodor Lynen fellowship of the
Alexander von Humboldt foundation.}
\begin{abstract}
  In \cite{BC1}, Buff and Ch{\'e}ritat proved that there are
  quadratic polynomials having Siegel discs with smooth
  boundaries. Based on a simplification of Avila, we give yet another
  simplification of their proof. The main tool used is a harmonic
  function introduced by Yoccoz whose boundary values are the sizes of
  the Siegel discs. The proof also applies to some other families of
  polynomials, entire and meromorphic functions.
\end{abstract}

\maketitle

\section{Introduction and Statement of Results}
In 1970 Katok asked whether there exist analytic circle
diffeomorphisms which are $\cC^\infty$ but not analytically
linearizable. This question was answered affirmatively by
P{\'e}rez-Marco in \cite{PM}, through the construction of Siegel discs
with smooth boundaries. P{\'e}rez-Marco also announced a proof of the
existence of smooth Siegel discs in the quadratic family $P_\lambda(z)
= \lambda z(1-z)$. Later Buff and Ch{\'e}ritat found a different proof
of this result \cite{BC1}, which was subsequently simplified by Avila
\cite{Av}, culminating in the joint paper \cite{ABC}. The present note
gives yet another simplification of the proof, based on representing
the size of the Siegel discs as boundary values of a certain harmonic
function in the unit disc, introduced by Yoccoz in \cite{Yo}. Yoccoz
used this function to give a short proof that Siegel discs exist for
almost all rotation numbers, and later Carleson and Jones used it to
prove that the critical point is on the boundary of the Siegel disc
for almost all rotation numbers. The method of Carleson and Jones was
later applied by Rippon in \cite{Ri} to show that almost all Siegel
discs in the exponential family $E_\lambda(z) = \lambda (e^z-1)$ are
unbounded.

The main theorem we are going to prove is the following.
\begin{thm}
  \label{thm:main}
  Let $f(z)=z+\ldots$ be meromorphic in $\C$ (not necessarily
  transcendental), and let $f_\lambda(z) = \lambda f(z)$. Assume that $f$ has only
  one non-zero critical or asymptotic value $v_0$. Then there is a
  dense subset $S \subset \R$ such that $f_{e^{2\pi i \alpha}}$ has a smooth Siegel
  disc for $\alpha \in S$.
\end{thm}
{\sc Remark.} In the case where $f$ has a rotational symmetry $f(\omega z)
= \omega f(z)$ with some root of unity $\omega$, it suffices to assume that $f$
has only one non-zero singular value up to multiplication by $\omega$. To
prove this if $\omega$ is an $n$-th root of unity, consider $F(w) =
f(w^{1/n})^n$. Then $F$ has at most one singular value, and $z \mapsto \lambda f(z)$
has a Siegel disc at 0 if and only if $w \mapsto \lambda^nF(w)$ has one.

We can easily check that the assumption is fulfilled for some families
and get the following immediate corollary.
\begin{cor}
  There are smooth Siegel discs in the following families:
  \begin{eqnarray*}
 P_{\lambda,d} (z) &=& \lambda \left( \left(1+\frac{z}{d} \right)^d - 1 \right) \\
 E_\lambda (z) &=& \lambda (e^z-1) \\
 f_\lambda (z) &=& \lambda z e^z \\
 g_\lambda (z) &=& \lambda \sin z \\
 h_\lambda (z) &=& \lambda \tan z
  \end{eqnarray*}
\end{cor}
Note that the polynomials $P_{\lambda,d}$ are conjugate to $z \mapsto z^d+c$ for
some $c=c(\lambda)$. 

Recently Buff and Ch\'eritat have refined the approximation technique
used to produce smooth Siegel disk in order to get Siegel disks with
any prescribed $C^\beta$ regularity\cite{BC2}.

\section{Proof}
We first recall some basic facts about linearization at attracting
fixed points, for a reference see e.g.~\cite{CG}.  For $0<|\lambda|<1$ the
function $f_\lambda$ has an attracting fixed point at 0, and by the
Koenigs-Poincar{\'e} linearization theorem there exists a linearizing map
$h_\lambda(z)=z+ O(z^2)$ with $h_\lambda (f_\lambda(z)) = \lambda h_\lambda (z)$ near 0. The
functional equation allows to extend $h_\lambda$ to the basin of attraction
of 0 for $f_\lambda$. We define the \emph{Yoccoz function} $w_\lambda$ of the
family $f_\lambda$ as $w(\lambda) = h_\lambda (\lambda v)$. Notice that $\lambda v$ is the unique
critical/asymptotic value of $f_\lambda$ with infinite forward orbit. The
inverse map $h_\lambda^{-1}$ of $h_\lambda$ is defined in the disc $\{ |w| < |w(\lambda)|
/ |\lambda| \}$ and does not extend to any larger disc centered at 0.  In
\cite{Yo}, Yoccoz proved that the radial limit $R(\alpha) =
\ds\lim_{r\to1}|w(re^{2\pi i \alpha})|$ exists for every $\alpha \in \R$. Furthermore, it
coincides with the conformal radius of the Siegel disc of $f_\lambda$ in the
case of linearizability, and it is 0 otherwise. (Yoccoz proved this
only for the case of the quadratic family, but his proof carries over
in our case without any essential changes.) Apart from basic facts
about conformal mappings, the only tool in Yoccoz's proof is the
Birkhoff Ergodic Theorem. By the Koebe 1/4-Theorem we get the estimate
$|w(\lambda)| < 4|v|$ for all $\lambda \in \D$, and the Koebe distortion theorem
also gives the asymptotic form $w(\lambda) = v \lambda + O(\lambda^2)$ at 0.
Furthermore, $w(\lambda)/ \lambda $ has no zeros in the unit disc, and it is still
bounded by $4|v|$. Thus $u(\lambda) = \log |w(\lambda)/ \lambda|$ is a harmonic function
in the unit disc, bounded above by $M:= \log 4 + \log|v|$. By Yoccoz's
result, the radial limit $\rho(\alpha) = \ds\lim_{r\to1} u(re^{2\pi i \alpha}) \in [-\infty, M] =
\log R(\alpha)$ exists for every $\alpha \in \R$.
\begin{prop}
  \label{prop:ivp}
  The function $\rho$ satisfies the following properties.
  \begin{itemize}
  \item [(a)] $\rho(\alpha) = -\infty$ for every $\alpha \in \Q$.
  \item [(b)] $\rho(\alpha)$ is finite for almost every $\alpha \in \R$.
  \item [(c)] $\ds \rho(\alpha) = \limsup_{\beta \searrow \alpha} \rho(\beta) = \limsup_{\beta \nearrow \alpha} \rho(\beta)$ for every $\alpha$.
  \item [(d)] $\rho$ satisfies the intermediate value property. 
\end{itemize}
\end{prop}

{\sc Proof.} Claim (a) immediately
follows from Yoccoz's result. Claim (b) follows from Fatou's Theorem
for positive harmonic functions. Claim (c) implies (d) by the
following argument: Let $\alpha<\beta$ and  $\rho(\alpha) < r < \rho(\beta)$. Define $\gamma_0 :=
\inf \{ \gamma \in [\alpha,\beta] : \rho(\gamma) \geq r \}$. Then there is a sequence $\gamma_n$ with
$\rho(\gamma_n) \geq r$ and $\gamma_n \to \gamma_0$, thus $\rho(\gamma_0) \geq r$ by (c). On the other
hand, $\ds \rho(\gamma_0) = \limsup_{\gamma\nearrow\gamma_0} \rho(\gamma) \leq r$, so $\rho(\gamma_0) = r$. The case
$\rho(\alpha)>r>\rho(\beta)$ is treated similarly.

Proof of (c): Upper semicontinuity of $\rho$ is clear. For every $\alpha \in \R$
with $\rho(\alpha) > -\infty$ there is a conformal map $g_\alpha(w) = w + \ldots$ defined in
$\{ |w| < R(\alpha) = e^{\rho(\alpha)} \}$ with $f_\lambda (g_\alpha (w)) = g_\alpha (\lambda w)$, where we write $\lambda =
e^{2\pi i \alpha}$. If $\alpha_n \to \alpha$ with $\rho(\alpha_n) \to r > -\infty$, then the sequence of
conformal maps $(g_{\alpha_n})$ is normal in $\{ |w| < e^r \}$ and every
subsequential limit $g(w) = w + \ldots$ is again a conformal map,
satisfying $f_\lambda (g(w)) = g(\lambda w)$, so $f_\lambda$ has a Siegel disc of
conformal radius $\geq e^r$, thus $\rho(\alpha) \geq r$.

The other direction is not quite as obvious, and the paper of Avila
uses some deep results of Risler to show it. However, it also follows
from the fact that $\rho$ can be represented as the boundary function of a
harmonic function in the unit disc. Fix $\alpha\in\R$ and let $\ds L :=
\limsup_{\beta \nearrow \alpha} \rho(\alpha)$ and $\ds R := \limsup_{\beta \searrow \alpha}$. For every $\epsilon > 0$
there exists $\delta > 0$ such that $\rho (\beta) \leq L+\epsilon$ for $\alpha - \delta \leq \beta < \alpha$ and
$\rho(\beta) < R+\epsilon$ for $\alpha < \beta \leq \alpha + \delta$.  Let $u_\epsilon$ denote the Poisson integral
of the function
\begin{equation*}
  u_\epsilon (e^{2\pi i \beta}) = \begin{cases}
    L+\epsilon & \text{ for $\alpha - \delta < \beta < \alpha$,} \\
    R+\epsilon & \text{ for $\alpha < \beta < \alpha + \delta$,} \\
    M & \text{ otherwise,}
  \end{cases}
\end{equation*}
in the unit disc. Then $u \leq u_\epsilon$ in the unit disc by the maximum
principle. This implies $\rho(\alpha) =\ds \lim_{r\to1} u(r e^{2\pi i \alpha}) \leq
\lim_{r\to1} u_\epsilon(r e^{2\pi i \alpha}) = \frac{L+R}2 + \epsilon$. As this holds for any
positive $\epsilon$, we get $\rho(\alpha) \leq \frac{L+R}2$. On the other hand upper semicontinuity
implies $\rho(\alpha) \geq \max(L,R)$. These two estimates together yield $\rho(\alpha) = L = R$. \qed

One immediate corollary is the following curious property, central to
the proof of the existence of smooth Siegel discs.
\begin{cor}
  \label{cor:curious}
  For any $\alpha\in \R$ and $ \tilde{\rho} < \rho(\alpha)$ there exists a sequence $\alpha_n \to \alpha$
  with $\rho(\alpha_n) = \tilde{\rho}$.
\end{cor}
{\sc Proof.} There exists a sequence of rational number $(\beta_n)$
converging to $\alpha$. We then have $\rho(\beta_n) = -\infty < \tilde{\rho} < \rho(\alpha)$, and
application of the intermediate value property yields the sequence
$\alpha_n$. \qed

The rest of the proof essentially follows Avila's presentation. Let
$F_r$ be the space of holomorphic functions in $\{ |w|<r \}$, endowed
with the compact-open topology. Define norms $\| f \|_r := \ds\sup_{k\geq
  0, |w|<r} \frac{|f^{(k)}(w)|}{[(k+2) \log (k+2)]^k}$ for holomorphic
functions in $\{ |w|<r \}$, and let $E_r$ be the space of functions for
which this is finite. Then $(E_r, \| . \|_r)$ is a Banach space of
holomorphic functions which are still $\cC^\infty$ in the closed disc $\{
|w| \leq r \}$. In fact, it is a class of quasi-analytic functions,
i.e.~if all derivatives of some $f \in E_r$ vanish at a point $w$ with
$|w|=r$, then $f\equiv0$. Furthermore, if $r>s$, then the inclusion mapping
$F_r \hookrightarrow E_s$ is continuous, as can easily be seen from the Cauchy
integral formula. 

Now let $\alpha_0\in\R$ be a number with $\rho_0=\rho(\alpha_0)>\infty$, let $\epsilon_0>0$ and $\rho_\infty
< \rho_0$ be arbitrary. We will find $\alpha_\infty \in (\alpha_0 - \epsilon_0, \alpha_0 + \epsilon_0)$,
$\rho(\alpha_\infty) = \rho_\infty$, and the corresponding Siegel disc of $f_{e^{2\pi i
    \alpha_\infty}}$ having smooth boundary.  Fix $\delta > 0$ such
that $\| g - g_{\alpha_0} \|_{r_\infty} \le 2\delta$ implies $g' \ne 0$
on the closed disk $\{ |w| \le r_\infty \}$. 
Let $(\rho_n)$ be a strictly
decreasing sequence converging to $\rho_\infty$. For convenience we will adopt
the notation $r_n = e^{\rho_n}$. We recursively define sequences $(\alpha_n)$
and $(\epsilon_n)$ as follows. Let $\alpha_{n+1} \in (\alpha_n - \epsilon_n, \alpha_n + \epsilon_n)$ with
$\rho(\alpha_{n+1}) = \rho_{n+1}$ and $\| g_{\alpha_n} -
g_{\alpha_{n+1}} \|_{r_\infty} \leq 2^{-n} \delta $.
This is possible because by Corollary \ref{cor:curious} there is a
sequence $(\alpha_{n,k})$ converging to $(\alpha_n)$ with $\rho(\alpha_{n,k}) =
\rho_{n+1}$. The corresponding linearizing maps $g_{\alpha_{n,k}}$ converge
locally uniformly on $\{ |w| < r_{n+1} \}$ to $g_{\alpha_n}$, so they
converge in $\| . \|_{r_\infty}$ because $r_\infty < r_{n+1}$. Thus we can pick
$\alpha_{n+1}$ among the $(\alpha_{n,k})$ to satisfy our requirements. Now
choose $\epsilon_{n+1}$ so that $[\alpha_{n+1} - \epsilon_{n+1}, \alpha_{n+1} + \epsilon_{n+1}] \subset
[\alpha_n - \epsilon_n, \alpha_n + \epsilon_n]$ and $\rho(\alpha) < \rho_n$ for $|\alpha - \alpha_{n+1}| \leq
\epsilon_{n+1}$. This is possible due to upper semicontinuity of $\rho$, implied
by Proposition \ref{prop:ivp}(c). Let $\alpha_\infty = \ds \lim_{n\to\infty} \alpha_n$ the
unique point in $\bigcap_n [\alpha_n - \epsilon_n, \alpha_n + \epsilon_n]$. Then by upper
semicontinuity, $\rho(\alpha_\infty) \geq \ds \lim_{n\to\infty} \rho(\alpha_n) = \rho_\infty$. On the other
hand, by construction we have that $\rho(\alpha_\infty) < \rho_n$ for all $n$, so
$\rho(\alpha_\infty) = \rho_\infty$. Furthermore, the sequence of linearizing maps
$(g_{\alpha_n})$ is a Cauchy sequence in $\| . \|_{r_\infty}$, so the limit
$g_{\alpha_\infty}$ is in $E_{r_\infty}$, and $\| g_{\alpha_\infty} -
g_{\alpha_0} \|_{r_\infty} \le 2 \delta$, so $g_{\alpha_\infty}' \ne
0$ on $\{ |w| \le r_\infty \}$. This implies that the Siegel disc of
$f_{e^{2\pi i \alpha_\infty}}$ is $S_\infty = g_\infty(\{ |w| < r_\infty \})$ and that it is a
domain with a $C^\infty$ boundary. \qed

\end{document}